\documentclass[11pt,english]{article}
\usepackage[T1]{fontenc}
\usepackage[latin9]{inputenc}
\usepackage[letterpaper]{geometry}
\usepackage{url}
\usepackage{amsmath}
\usepackage{amsthm}
\usepackage{amssymb}
\usepackage{bm}

\makeatletter

\usepackage{dsfont}
\usepackage{babel}
\usepackage{bbm}

\numberwithin{equation}{section}
\theoremstyle{plain}
\newtheorem*{thm*}{Theorem}
\newtheorem*{lem*}{Lemma}
\newtheorem*{prop*}{Proposition}

\makeatother

\newcommand {\N}        {\mathbb{N}}
\newcommand {\C}        {\mathbb{C}}
\newcommand {\R}        {\mathbb{R}}
\newcommand {\ind}[1] 	{\mathds{1}\{#1\}}
\newcommand {\ip}[1]  	{\langle#1\rangle}
\newcommand {\E}        {\mathbb{E}}
\newcommand {\cb}       {\boldsymbol{c}}
\newcommand {\n}        {\boldsymbol{n}}
\newcommand {\s}        {\boldsymbol{s}}
\newcommand {\bS}       {\boldsymbol{S}}
\newcommand {\tb}       {\boldsymbol{t}}

\newcommand {\z}        {\boldsymbol{z}}
\newcommand {\vloop}  	{\bm \omega}
\newcommand {\ab}       {\boldsymbol{a}}
\newcommand {\tbold}    {\boldsymbol{t}}
\newcommand {\btheta}	  {\bm \theta}
\newcommand {\pinV}[1]	{\prod_{#1 \in V}}
\newcommand {\paths}[1]	{\mathcal{P}_#1}
\newcommand {\curs}[1]	{\mathcal{C}_#1}
\newcommand {\eL}		    {\mathcal{L}}
\newcommand {\fnset}[1]	{\N_{\text{fin}}^#1}
\newcommand {\size}[1]	{\left| #1 \right|}

\begin{document}

\title{Weighted graphs and complex Gaussian free fields}
\author{Gregory F. Lawler
	\thanks{Research supported by NSF grant DMS-1513036} \\
		Petr Panov}
\maketitle

\begin{abstract}
We prove a combinatorial lemma about the distribution of directed currents
in a complex ``loop soup'' and use it to give a new proof of the isomorphism
relating loop measures and complex Gaussian fields.
\end{abstract}

\section{Introduction}
Loops and related measures are useful tools in the analysis of random walks.
They have recently come under study in \cite{key5} as a discrete analogue of
the Brownian loop soup introduced in \cite{key4}, which itself was motivated by
the study of the Schramm-Loewner evolution. Such measures were also explored in
a continuous setting by Yves Le Jan \cite{key2}. One of his findings was the
connection between the Gaussian free field and the occupation field of a
Poissonian ensemble of Markov loops. This connection can be viewed as a version
of the Dynkin's isomorphism theorem \cite{key6}. In \cite{key7} this isomorphism
was extended to connect certain non symmetric Markov processes and complex
Gaussian fields. 

A version of the isomorphism theorem using the discrete time loop soup was 
proved in \cite{key1} and \cite{key3}. Random walk on a finite graph can be 
fully described by a  substochastic transition matrix $Q$. Any event in this 
setting is essentially a union of chain trajectories, and its probability is an 
additive function on sets of trajectories, which is related to $Q$. It is not 
uncommon in statistical physics to interpret events and their probabilities as
configuration collections and weights, respectively. Even if $Q$ takes complex
values, in some cases we can still build meaningful objects that have
probabilistic analogues, such as loop soups, by putting potentially complex
weights on paths. That proof involves the loop soup at intensity $1/2$ and uses 
undirected currents.

This generalization is one of the main points of the discussion in \cite{key3}.
The complex Gaussian free field is introduced there as a pair of real Gaussian
free fields with potentially negative correlations between fields and within
each field. A version of the isomorphism theorem is formulated and proved there
by comparing the Laplace transforms of a complex Gaussian field squared and a
continuous occupation field of a complex loop soup. A combinatorial proof of
the isomorphism can be found in \cite{key1}, which is one the one hand
discussed under the assumption that weights correspond to a certain probability
space, and on the other hand does not use that assumption in a significant way.

This note heavily relies on, and serves as a continuation those two papers. Here
we adapt the arguments from \cite{key1} to the complex setting and extend some
of the results of \cite{key3} to a wider range of weights. The key new results
here are the exact distribution on directed currents, induced by the random 
walk loop soup at intensity 1 (which is presented in our Proposition), and the 
isomorphism theorem, which connects the continuous occupation field of the loop 
soup and the absolute value of a complex Gaussian free field squared. To prove 
the latter, we do not utilize the Laplace transform, which is commonly used to 
show the isomorphism in the literature. Another advantage of the theorem 
presented here is that it involves the random walk loop soup measure at 
intensity $1$, which is easier to analyze than the loop soup at intensity 
$1/2$. Hopefully, our proof sheds some light on a seemingly accidental 
connection between the loop measures and the Gaussian free field. The 
isomorphism theorems proved in \cite{key1} and \cite{key3} become a special 
case of our Theorem.

This paper is structured as follows. We first introduce the setup and basic
notations. Then we state the main results of the paper, including the
isomorphism theorem. All the proofs are contained in the final section.
\section{Basic definitions}
Consider a finite complete digraph $G = (V, E)$ with $N = |V|$ vertices. 
Directed edges $E \cong V \times V$ are identified with ordered pairs of
vertices; note that we allow self-edges. The set of vertices is ordered:
$V = (v_j)_{j=1}^N := (v_1, v_2, \ldots, v_N)$ and whenever we take an ordered 
subset of it, we preserve the order between the vertices. We will often 
use the following ordered subsets: $V_k = (v_j)_{j=k}^N$, for $k \in [N]$. For 
future reference, we use  $(\cdot)$ and  $\{\cdot\}$ to denote ordered and 
unordered sets, respectively.

We call functions $q:E \to \C$ \emph{weights} on directed edges. For any
$u, v \in V$, we write $q_{uv}$ instead of $q(u,v)$ for brevity. Let
$Q = (q_{uv} : u, v \in V) \in \C^{N \times N}$, and call $q$ a
\emph{Hermitian weight}, if $Q$ is Hermitian. We say that $q$ is 
\emph{integrable}, if $\rho(|Q|) < 1$, where $\rho$ denotes the spectral radius 
and $|Q| = \bigl( |q_{uv}| : u, v \in V \bigr)$. The notation $|Q|$ always 
refers to the matrix of absolute values; to denote the matrix determinant, we 
use $\det$. If $U \subseteq V$, then $Q_U$ denotes the restriction of $Q$ to 
rows and columns that correspond to vertices in $U$, that is, $Q_U = (q_{uv} : 
u, v \in U) \in \C^{|U| \times |U|}$.

A \emph{path} $\omega$ of \emph{length} $\size{\omega} = k$ in $U \subseteq V$ 
is a sequence of $k + 1$ vertices in $U$:
\[
  \omega = (\omega^0, \omega^1, \ldots, \omega^k) =
  (\omega^j)_{j=0}^k, \qquad \{\omega^j\}_{j=0}^k \subseteq U.
\]
Paths of length $0$ are called \emph{trivial}. Equivalently, a path $\omega$ is 
a sequence of $k = \size{\omega}$ directed edges, where the second vertex of
each edge is the same as the first vertex of the next edge in the sequence:
\[
  \omega = \oplus_{j=1}^k e_j, \qquad e_j = (\omega^{j-1},\omega^j) \in E \qquad
  \forall j \in [k] := \{1, 2, \ldots, k\}.
\]
Concatenation of edges is represented by $\oplus$ here. It can be applied to
paths similarly, as long as each next path in the sequence that we concatenate
starts where the previous path has ended.

If $\{u, v\} \subseteq U$, let $\paths{U}(u, v)$ denote the set of paths in $U$
starting at $u$ and ending at $v$. Paths $\paths{U}(v) := \paths{U}(v, v)$ are 
called \emph{loops rooted at} $v$ in $U$ and contain the \emph{trivial loop}
consisting of a single vertex. We use $\paths{U}$ to denote
$\cup_{u, v \in U}\paths{U}(u, v)$. We have previously defined $q$ as a 
function on directed edges, and we will also use it to denote the following 
function on $\paths{V}$:
\[
  q(\omega) = \prod_{j=1}^k q(\omega^{j-1}, \omega^j) =
  \prod_{j=1}^k q(e_j) \quad \text{for} \quad
  \omega = (\omega^j)_{j=0}^k = \oplus_{j=1}^k e_j \in \paths{V}, \quad
  k = \size{\omega} > 0,
\]
and we let $q$ be equal to $1$ on the trivial loops in $\paths{V}$.  Note that 
if $q$ is integrable, then it defines a complex measure on $\paths{V}$.
 
For any path $\omega \in \paths{V}$, we let
$\n(\omega) = \{n_u(\omega): u \in V\}$ denote the {\em (vertex) local time},
where
\[
  n_u(\omega) = \sum_{j=1}^{\size{\omega}} \delta_u(\omega^j):=
  \sum_{j=1}^{\size{\omega}} \ind{\omega^j = u}.
\]
Note that we are not counting the visit at time $0$. In particular,
$\n = \boldsymbol{0}$ on trivial loops. We write
$\cb = \bigl(c_{uv}(\omega) : u, v \in V \bigr)$ for the
{\em (directed) edge local time}, where
\[
  c_{uv}(\omega) = \sum_{j=1}^{\size{\omega}}
  \ind{\omega^{j-1} = u, \omega^j = v}.
\]
A matrix $(C_{uv})_{u, v \in V}$ with entries in $\N = \{0, 1, 2, \ldots\}$ is
called a {\em (directed) current}, if
\[
\sum_{v \in V} C_{uv} = \sum_{v \in V} C_{vu}, \qquad \forall u \in V.
\]
If $U\subseteq V$, we use $\curs{U}$ to denote the set of currents restricted
to $U$, that is, such that $C_{uv} = 0$ if either $u$ or $v$ is in
$V \setminus U$. Note that when $\omega$ is a rooted loop, the matrix
$C = \cb(\omega)$ is a current. It is then immediate, that for 
each $u \in V$,
\begin{equation} \label{ndef}
  n_u(\omega) = \frac12 \, \sum_{v \in V}
  \bigl[ c_{uv}(\omega) + c_{vu}(\omega) \bigr].
\end{equation}
We can extend the definition of $q$ to currents. If $C \in \curs{V}$, we set
\[
  q(C) = \pinV{u, v} q_{uv}^{C_{uv}}.
\]
An {\em (oriented) unrooted loop} is an equivalence class of nontrivial rooted 
loops under cyclic permutations:
\[
  (\omega^0, \omega^1, \ldots, \omega^k = \omega^0) \sim
  (\omega^1, \omega^2, \ldots \omega^{k-1}, \omega^0, \omega^1) \sim \cdots \sim
  (\omega^{k-1}, \omega^0, \omega^1, \ldots, \omega^{k-2}, \omega^{k-1}).
\]
The set of unrooted loops is denoted by $\eL$. If a rooted loop $\omega$
represents $l \in \eL$, we will write $\omega \in l$. The set of unrooted loops 
whose representatives stay in $U \subseteq V$ and visit $v \in U$ at least once 
is denoted by $\eL_U(v)$.  The definitions of $q$, $\n$ and $\cb$ are extended
from $\paths{V}$ to $\eL$ by taking any rooted representatives:
\[
  q(l) = q(\omega),\ \n(l) = \n(\omega),\ \cb(l) = \cb(\omega), \quad
  \text{ if } \quad \omega \in l,\quad \forall l \in \eL.
\]
Such an extension does not depend on the choice of the representative.

If $\mathcal{X}$ is any countable set, we let $\fnset{\mathcal{X}}$ stand for
finite multisets of elements from $\mathcal{X}$, that is, the set of functions
$\mathcal{X} \mapsto \N$, which are supported on a finite set. Local times $\n$
and currents $\cb$ can be viewed as functions on $\fnset{\eL}$:
\[
  \n(\s) = \sum_{l \in \eL} \n(l) \cdot s_l, \qquad 
  \cb(\s)= \sum_{l \in \eL} \cb(l) \cdot s_l, \qquad \s \in \fnset{\eL}
\]
\section{Main results}
\subsection{Loop measures and occupation fields}
If $q$ is an integrable  weight on $V$, we define the
\emph{unrooted loop measure} $m$ by
\begin{equation} \label{mdef}
  m(l) = \sum_{\omega \in l} \frac{q(\omega)}{\size{\omega}} = 
  \frac{q(l)}{d(l)},
\end{equation}
where $d(l)$ is the largest integer $d$ such that every representative of $l$ 
consists of the concatenation of $d$ identical loops. If $q$ is integrable, $m$
is a complex measure on $\eL$.

The \emph{(random walk) loop soup (at intensity 1)} is a collection of
independent Poisson random variables indexed by $\eL$ with intensity
$e^{-m(l)}$. If $q$ is complex, we interpret this as the measure on 
finite multisets of unrooted loops:
\begin{equation} \label{numdef}
  \nu_m\{\s\} = \prod_{l \in \eL} \frac{e^{-m(l)} \, m(l)^{s_l}}{s_l!} =
  e^{-m(\eL)} \prod_{l \in \eL} \frac{m(l)^{s_l}}{s_l!},\qquad
  \forall \s \in \fnset{\eL}.
\end{equation}
We write $\nu_c$ and $\nu_*$ for the pushforwards of this measure as measures
on $\curs{V}$ and $\N^V$:
\[
  \nu_c\{C\} = \sum_{\s \rightarrow C} \nu_m\{\s\}, \qquad
  \nu_*\{\n '\} = \sum_{\s \rightarrow \n '} \nu_m\{\s\},
\]
where the sums are over all $\s \in \fnset{\eL}$ that produce the current $C$ 
and the local time $\n '$, respectively. We call $\nu_c$ the
\emph{(directed) current field} and $\nu_*$ the
\emph{discrete occupation field}.

Our first result, gives the distribution for the current field of any 
integrable weight (not necessarily Hermitian).

\begin{prop*}
  If $q$ is an integrable weight and  $C \in \curs{V}$, then
  \begin{equation} \label{curdist}
    \nu_c(C) = \det(I - Q) \, q(C) \pinV{u} \frac{n_u(C)!}{\pinV{v} C_{uv}!}.
  \end{equation}
\end{prop*}

The proof of this fact is combinatorial in nature and revolves around the 
identity \eqref{comb}, which can be viewed as a useful result on its own.

Given a discrete occupation field, the \emph{continuous occupation field} is
obtained by independently at each vertex $u$ replacing $n_u$ with the sum of
$n_u + 1$ independent exponential random variables with mean one. We write this
distribution as $\nu_n$.  We can give its density with respect to Lebesgue 
measure $\lambda_N$ on $\R_+^N$:
\[
  \frac{d\nu_n}{d\lambda_N}(\tbold) = \sum_{\s \in\fnset{\eL}}
  \left[ \nu_m\{\s\} \pinV{u} \frac{t_u^{n_u(\s)}e^{-t_u}}{n_u(\s)!} \right],
  \qquad \tbold \in \R_+^N
\]

Due to the proposition above, this can be written as
\begin{eqnarray}
  \frac{d\nu_n}{d\lambda_N}(\tbold) & = &
  \sum_{C \in \curs{V}} \left[ \nu_c(C) \, \pinV{u}
  \frac{t_u^{n_u(C)} e^{-t_u}}{n_u(C)!} \right]  \nonumber \\  & = &
  \det(I - Q) \sum_{C \in \curs{V}} \left[ q(C) \, \pinV{u}
  \frac{t_u^{n_u(C)} e^{-t_u}}{\pinV{v}C_{uv}!} \right] \label{nundef}
\end{eqnarray} 
\subsection{Bubble soup}
In order to prove the Proposition, we define and analyze certain 
auxiliary measures. For $v \in U \subseteq V$, a \emph{growing loop} in $U$ at 
$v$ 
induced by $q$
(at time $t = 1$) is a ``random'' rooted loop in $\paths{U}(v)$
sampled as follows. If $\nu_g$ denotes the measure on the growing loop on 
$\paths{U}(v)$, then
\[
  \nu_g(\omega) = \frac{q(\omega)}{G_U(v,v)}, 
\]
where $G_U(v, v)$ denotes the Green's function. Note that
\begin{equation} \label{expm}
  G_U(v,v) = \sum_{\omega \in \paths{U}(v)} q(\omega) =
  \exp \left\{ \sum_{l \in \eL_U(v)} m(l) \right\}.
\end{equation}
Indeed, the first equation follows from a standard renewal argument and uses 
the fact that $q$ is a complex measure on paths with finite total variation. 
For the final expression, see Lemma 3.1 in \cite{key3}.

The {\em bubble measure} $\nu_b$ is the measure  on  $N$-tuples
$\vloop = (\omega_j: \omega_j \in \paths{{V_j}}(v_j))_{j=1}^N$
given by the product measure
\begin{equation} \label{bubdef}
  \nu_b(\vloop) = \prod_{j=1}^N \nu_g(\omega_j) =
  \prod_{j=1}^N \frac{q(\omega_j)}{G_{V_j}(v_j, v_j)} =
  \frac{q(\vloop)}{\det G} , \quad \text{where} \quad q(\vloop) =
  \prod_{j=1}^N  q(\omega_j).
\end{equation}
Here we have used the following well-known formula (for example, see the
Proposition 3.5 in \cite{key1}):
\begin{equation} \label{detG}
  \det G = \prod_{j=1}^N G_{V_j}(v_j, v_j).
\end{equation}
Note that the definition depends on the ordering of $V$, but if we forget the 
order in $\vloop$, then it is immediate from $\eqref{bubdef}$, that the 
resulting measure will not depend on the order in $V$.

The following statement allows us to work with bubble soup instead of the 
unrooted loop soup to derive the current distribution \eqref{curdist}.

\begin{lem*}
	For any ordering of $V$, the measure induced on currents by $\nu_b$ is
	$\nu_c$.
\end{lem*}

This follows immediately from the Proposition 5.8 of \cite{key1} for general
intensities when $Q$ is a substochastic matrix. In the case of intensity one 
and positive weights a similar result was established in the Proposition 9.4.1 
of \cite{key8}. Unfortunately, there was a misstatement in the latter proof of 
Problem 9.1 which was part of the proof. Because of this unfortunate misprint, 
we will redo the proof here. We will also show that the argument applies to 
general integrable weights.
\subsection{Isomorphism theorem}
We now assume that $q$ is Hermitian, and thus $G = (I - Q)^{-1}$ is a positive
definite Hermitian matrix. The
{\em (discrete centered) complex Gaussian free field} $Z = (Z_v : v \in V)$ on 
$V$ with covariance $G$ is a random complex vector in $\C^N$ with density
\[
  f_Z(\z) = \frac{\exp \bigl\{ -\ip{\z,G^{-1}\z} \bigr\} }{\pi^N\det G}, \qquad
  \z \in \C^N
\]
with respect to the Lebesgue measure on $\C^N$; here $\ip{\cdot, \cdot}$
denotes the dot product of complex vectors. $Z$ satisfies the following
covariance relations:
\[
  \E[\bar{Z}_u Z_v] = G(u, v), \qquad \E[Z_u Z_v] = 0; \qquad
  \forall u, v \in V.
\]
We can decompose the Green's function into the real and imaginary parts:
$G = G^R + i G^I$. Since $G$ is Hermitian, $G^R$ is symmetric and $G^I$ is
antisymmetric. A complex Gaussian free field on a set of $N$ elements can be
viewed as a real field on $2N$ elements. Indeed, let
\begin{equation} \label{gffreal}
  (Z', Z'') := (Z_u')_{u \in V} \oplus (Z_u'')_{u \in V} \sim \mathcal{N}
  \left( 0, \Big(
  \begin{array}{cc}
    G^R & -G^I\\
    G^I & G^R
  \end{array} \Big) \right),
\end{equation}
where $\oplus$ denotes the concatenation of sequences. According to the
Proposition 4.5 in \cite{key3},
\begin{equation} \label{zeid}
  Z \overset{\mathcal{D}}{=} (Z' + i Z'') / \sqrt{2},
\end{equation}
that is, the probability distributions of these complex random vectors are the
same.

Let $f_{|Z|^2}$ denote the density of $|Z|^2 = (\bar{Z}_u Z_u)_{u \in V}$ with
respect to Lebesgue measure $\lambda_N$ on $\R_+^N$. According to the
Theorem 2 in \cite{key3}, we should expect that the continuous occupation field
at intensity $1$ has the same density as the square of the absolute value of a
complex Gaussian free field. The following generalizes the isomorphism theorems
as stated in \cite{key1} and \cite{key3}.

\begin{thm*}
  If $q$ is an integrable, Hermitian weight, then the continuous occupation 
  field $\nu_n$ has the same distribution as $|Z|^2/2$ where $Z$ is a complex
  Gaussian free field with covariance matrix $G = (I - Q)^{-1}$.
\end{thm*}

In view of \eqref{zeid}, this result can be interpreted differently. If $Z'$ and
$Z''$ are two real Gaussian free fields with correlation structure as in
\eqref{gffreal}, then $\nu_n$ has the same distribution as
$\bigl( |Z'|^2 + |Z''|^2 \bigr) /2$.

If $Q$ is a nonnegative integrable weight, then the distribution of $|Z|^2$ is
the same as that of $|X|^2 + |Y|^2$ where $X, Y$ are independent real Gaussian
fields with covariance matrix $(I - Q)^{-1}$. In this case, the result above 
reduces to the usual isomorphism theorem, which states that $|X|^2/2$ has the
same distribution  as the continuous occupation field at time $1/2$.
\section{Proofs}
\subsection{Proof of the Lemma}
We start by combining \eqref{expm} and \eqref{detG} to see that
$\det G = e^{m(\eL)}$. In view of that, and also \eqref{numdef} and
\eqref{bubdef}, we see that the goal is to prove that
\begin{equation} \label{l1}
  \sum_{\vloop \rightarrow C} q(\vloop) =
  \sum_{\s \rightarrow C} \left[ \prod_{l \in \eL}
  \frac{m(l)^{s_l}}{s_l!} \right],
\end{equation}
where the first sum is over all such tuples $\vloop = (\omega_j)_{j=1}^N$, that
$\sum_{j=1}^N \cb(\omega_j) = C$.

Let $\eL_k = \eL_{V_k}(v_k)$ for every $k \in [N]$, for brevity. Since 
$\{\eL_j\}_{j=1}^N$ are disjoint, any multiset $\s \in \fnset{\eL}$ can be 
uniquely represented
by multisets $\{\s^j\}_{j=1}^N$, where $\s^j \in \fnset{{\eL_j}}$ for each $j 
\in [N]$. Using \eqref{mdef}, we rewrite \eqref{l1} as
\begin{equation} \label{l2}
  \sum_{\vloop \rightarrow C} \left[ \prod_{j=1}^N q(\omega_j) \right] =
  \sum_{\s \rightarrow C} \left[ \prod_{j=1}^N \prod_{l \in \eL_j}
  \frac{q(l)^{s_l}}{s_l! \, d(l)^{s_l}} \right],
\end{equation}
Fix any $j \in [N]$. Let $u = v_j$, $L = \eL_j$ and $P = \paths{{V_j}}(u)$ for 
brevity. We shall now construct a mapping from $\fnset{L}$ to $P$. Take any
$\s \in \fnset{L}$ and order the unrooted loops in it arbitrarily. For each 
unordered loop, choose a representative loop in $P$, choosing uniformly at 
random from all the possibilities. Concatenate all the rooted loops in the 
order they were produced into a rooted loop $\omega$. If $o$ is the combination 
of ordering and choice of rooted loops, then we define $\psi(\s, o)$ to be the 
resulting loop $\omega \in P$. Let  $O(\s)$  denote the set of all the possible 
choices $o$ for the multiset $\s$. According to this definition,
\begin{equation}  \label{l3}
  \size{O(\s)} = \frac{S_{\s}!}{\prod_{l \in L} s_l!} \cdot
  \prod_{l \in L} \left( \frac{n_u(l)}{d(l)} \right)^{s_l} =
  S_{\s}!  \cdot \prod_{l \in L}  \frac{n_u(l)^{s_l}}{s_l! \, d(l)^{s_l}},
\end{equation}
where, as before, $d(l)$ is the largest integer such that any rooted 
representative of $l$ is a concatenation of $d(l)$ identical rooted loops, and 
$S_{\s} = \sum_{l \in L}s_l$.

Note that $\prod_{l \in L}q(l)^{s_l} = q(\omega)$ whenever 
$\psi(\s, o) = \omega$ for some $o \in O(\s)$.  We can now see from \eqref{l2}
and \eqref{l3}, that it is sufficient to prove that for any $\omega \in P$ with
$n_0 = n_u(\omega) \geq 1$,
\begin{equation} \label{cycle}
  \sum \frac{1}{S_{\s}! \prod_{l \in L}n_u(l)^{s_l}} = 1,
\end{equation}
where the sum is over all pairs $(\s, o)$ with $o \in O(\s)$ and 
$\psi(\s, o) = \omega$.

There is a natural bijection between $\psi$ and finite sequences of positive 
integers $(n_j)_{j=1}^k$ with $\sum_{j=1}^k n_j = n_0$, which we call seq$(k, 
n_0)$. Multiplying both sides of \eqref{cycle} by $n_0!$, we see that it is 
equivalent to the identity
\[
  \sum_{k=1}^\infty \sum_{\text{seq}(k, n_0)} \frac{n_0!}{k! \,
  \prod_{j=1}^k n_j} = n_0!
\]
To establish this we need to show that the left-hand side equals the number of 
permutations of $n_0$ elements. To see this, suppose $(n_j)_{j=1}^k$ are given 
and $(a_j)_{j=1}^{n_0} = (a_1, a_2, \ldots, a_{n_0})$ is a permutation of
$(1, 2, \ldots, n_0)$. Then we get another permutation by putting parentheses 
down:
\[
  (a_1, \ldots, a_{n_1}) ,\, (a_{n_1 + 1}, \ldots, a_{n_1 + n_2}) ,\, 
  \ldots, \, (a_{n_1 + n_2 + \ldots + n_{k-1} + 1}, \ldots, a_{n_0} ),
\]
and viewing this as a representation of a permutation by its cycle structure.
However, there are many ways to represent the same permutation. There are $k!$ 
ways to permute the elements of $(n_j)_{j=1}^k$, and for the cycle 
corresponding to $n_j$ there are $n_j$ choices for which element to call 
$a_{n_1 + n_2 + \ldots + n_{j - 1} + 1}.$ This establishes our claim.
\subsection{Proof of the Proposition}
We will prove this by induction on the number of vertices $N = \size{V}$, 
viewing the current measure as the pushforward of a bubble soup under the 
mapping $\cb$. If $C_{uv} \neq 0$ for some $u, v$ with $q_{uv} = 0$, then both 
sides of \eqref{curdist} equal zero. Hence we will assume that $C$ is a current 
such that $q_{uv} \neq 0$ if $C_{uv} \neq 0$. 

If $V = \{x\}$ is a singleton with $q = q_{xx}$, then the bubble soup consists
only of self-loops $\omega_k = (x, x, \ldots, x)$ at $x$ in $V$ with
$\size{\omega_k} = k$ and 
\[
  \nu_g\{\omega_k\} = (1 - q) \, q^k = \det(I - Q) \, q^k
\]
and $k \in \N$. A current $C$ with $C_{xx} = k$ satisfies $n_x(C) = k$ and only 
$\omega_k$ can induce it. Hence,
\[
  \nu_c\{C\} = \nu_b\{\omega_k\} = \det(I-Q) \, q^k.
\]
Therefore \eqref{curdist} holds in this case.
 
Now suppose that $V$ has $N = \size{V} \geq 2$ vertices, $x \in V$, and let
$U = V\setminus\{x\}$. The induction assumption is that \eqref{curdist} holds 
for the currents $C \in \curs{U}$, where $Q_U$ denotes $Q$ restricted to $U$. 
We call this measure $\nu_0$. Using \eqref{detG}, we see that 
\begin{equation} \label{p1}
  \det(I - Q) = \det(I - Q_U) / G_V(x, x)
\end{equation}
If we order the vertices in $V$ in such a way that $x$ is the first vertex,
then the construction of a bubble soup in $V$ starts by growing loops at $x$ in
$V$, and the growing loops that follow are fully contained in $U$. Let $\nu_+$
be the measure on currents induced by a growing loop at $x$ in $V$:
\begin{equation} \label{p2}
  \nu_+\{C^+\} = \sum_{\omega \in L(C^+)} q(\omega) / G_V (x, x),
\end{equation}
where $L(C^+)$ is the set  of loops in $\paths{V}(x, x)$, that induce a current
$C^+$:
\[
  L(C^+) = \bigl\{ \omega \in \paths{V}(x, x) :\cb(\omega) = C^+ \bigr\}.
\]
Note that for any $\omega \in \paths{V}$ we have
$q(\omega) = q\bigl(\cb(\omega)\bigr)$, therefore all the summands in \eqref{p2}
are equal, since they correspond to the same current. We can now rewrite 
\eqref{p2} in a simplified form:
\begin{equation} \label{p3}
  \nu_+\{C^+\}= W(C^+) \, q(C^+) / G_V(x, x),
\end{equation}
where $W(C^+) =  \size{L(C^+)}$. To get the distribution on currents induced by 
a bubble soup in $V$, we calculate the measure $\nu_{0}$ induced by a bubble 
soup in $U$, the measure $\nu_+$ induced by a loop growing at $x$ in $V$, and 
take their convolution:
\begin{equation}  \label{p4}
  \nu_c\{C\} = (\nu_+ \ast \nu_0)\{C\} = \sum_{(C^+, C^0) \in P_C}
  \nu_+\{C^+\} \cdot \nu_0\{C^0\},
\end{equation}
where $C \in \curs{V}$ and
\[
  P_C = \Big\{ (C^+, C^0): C^+ \in \curs{V}, C^0 \in \curs{U},
  C^+ + C^0 = C \Big\}.
\]
Since $q(C^+ + C^0) = q(C^+) \, q(C^0),$  we can combine \eqref{curdist},
\eqref{p3} and \eqref{p4} to see that it suffices to prove the
following combinatorial statement:
\[
  \pinV{u} \frac{n_u(C)!}{\pinV{v} C_{uv}!} = \sum_{(C^+, C^0)\in P_C}
  \left( W(C^+) \, \prod_{u \in U}\frac{n_u(C^0)!}
  {\prod_{v \in U}C_{uv}^0!} \right), \qquad \forall C \in \curs{V}
\]
Note that the products can be written as multinomial coefficients:
\begin{equation} \label{comb}
  \pinV {u} \binom{n_u(C)}{\{C_{uv}\}_{v \in V}} = \sum_{(C^+, C^0) \in P_C}  
  \left[ W(C^+) \, \prod_{u \in U} \binom{n_u(C^0)}{\{C_{uv}^0\}_{v \in U}}
  \right]
\end{equation}
Let us fix an ordering of $V$ starting with $x$ and we use the same ordering on 
$U$ (ignoring $x$). For every $u \in V$, let $n^u = n_u(C)$ and let
$S^u = S^u(C)$ be the set of $n^u$-tuples $\ab^u = (a_1^u, \ldots, a_{n^u}^u)$
in $V^{n^u}$ that contain $C_{uv}$ elements $v$ for every $v \in V$. Let
$\bS(C) = \bigl( S^u \bigr)_{u\in V}$ be the collection of such sequences. Note
that the left-hand side of \eqref{comb} is equal to
$\size{\bS} = \pinV{u} \size{S^u}$.

Now take any pair $(C^+, C^0) \in P_C$ and let $n_+^u$, $S_u^+$, $n_0^u$ and 
$S^u_0$ be the corresponding quantities. Define
\[
  \bS ' = \bS '(C) = \bigcup_{(C^+, C^0) \in P_C} \left[ L(C^+) \times \bS (C^0)
  \right].
\]
The right-hand side  of \eqref{comb} is $\size{\bS '}$, thus it suffices to 
give a bijection between $\bS$ and $\bS '$.
    
Suppose $(\ab^u)_{u \in V} \in \bS(C)$ are given. To map $\bS$ to 
$\bS '$, we define $\omega \in \paths{V}(x, x)$ by means of an algorithm.
\begin{itemize}
  \item Set $\omega = (x)$.  If $n^x = 0$, stop and output the trivial loop.
  
  \item  Otherwise, let $\omega = (x, a_1^x)$, remove $a_1^x$ from $\ab^x$ and 
  reset $n^x \rightarrow n^x - 1$. 
\end{itemize}  
For $j = 1, 2, \ldots$, we do the following.
\begin{itemize}
  \item If $\omega^j = x$ and $n^x = 0$, stop and output 
  $\omega = (\omega^0, \ldots, \omega^j)$ and $\{\ab^u : u \in U\}$.

  \item Otherwise, if $\omega^j = u$, let $\omega^{j + 1}$ equal $a_1^u$, 
  remove $a_1^u$ from $\ab^u$, and  reset $n^u \rightarrow n^u -1$.
\end{itemize}
If the algorithm is correct, then clearly $\omega \in L(C^+)$ for
$C^+ = \cb(\omega)$, also $C^0 := C - C^+ \in \curs{U}$ and
$(\ab^u)_{u \in V} \in \bS(C^0)$. The correctness follows from the current 
property of $C \in \curs{V}$. We cannot encounter a situation where
$\omega^j = u \neq x$ and $n^u = 0$, because that would imply that
\[
  \sum_{v \in V} C_{uv} = \size{\ab^u} < \sum_{v \in V}
  \sum_{k=1}^{\size{\ab^v}} \delta_u(a_k^v) = \sum_{v \in V} C_{vu}
\]
To get the inverse mapping, run the algorithm in reverse, that is, map an
$\omega \in L(C^+)$ to an element in $\bS (C^+)$ and concatenate it with a 
vector from $\bS (C^0)$.
\subsection{Proof of the Theorem}
Fix an ordering $V$. To avoid cumbersome notation, identify vertices with 
integers: $V = (1, 2, \ldots, N)$.  For $\z \in \C^N$, let $x = \text{Re}[\z]$ 
and $y = \text{Im}[\z]$. If we do the change of variables
$x_j + i y_j = \sqrt{t_j}\, e^{i \theta_j}$,
$t_j \in \R_+$ and $\theta_j \in [0, 2\pi)$, we get
\[
  f_Z(\tb, \btheta) = g(\tb) \exp \left\{ \sum_{j,k=1}^N \sqrt{t_j t_k} \,
  q_{jk}\, e^{i(\theta_k - \theta_j)}\right\} ,
\]
where
\[
  g(\tb) := \exp \left\{ -\sum_{j=1}^N t_j \right\} \frac{\det(I - Q)}{(2\pi)^N}
\]
To get the marginal density $f_{|Z|^2}(\tb)$ for $\tb \in \R_+^N$, we integrate 
over $\btheta\in T:= [0, 2\pi)^N$: 
\begin{equation} \label{gffdens} 
  f_{|Z|^2}(\tb) = g(\tb) \int_T \exp \left\{ \sum_{j,k=1}^N \sqrt{t_j t_k}  \,
  q_{jk} \, e^{i(\theta_k - \theta_j)}\right\} \, d \btheta
\end{equation}
Next we find the density of the occupation field using its current
representation. Suppose that we have a matrix $C \in \N^{N \times N}$. Then
\[
  C \in \curs{V} \iff \sum_{j=1}^N (C_{jk} - C_{kj}) = 0 \quad \forall k \in [N]
\]
Since the right-hand side is always an integer for $C \in \N^{N \times N}$, we
see that  
\begin{eqnarray}
  \ind{C \in \curs{V}} & = & \prod_{j=1}^N \int_0^{2\pi} \frac{d\theta_j}{2\pi}
  \exp \left\{ i \theta_j \sum_{k=1}^N (C_{kj} - C_{jk})\right\} \nonumber \\
  & = & \int_T \frac{d\btheta}{(2\pi)^N}
  \left[ \prod_{j,k=1}^N e^{i C_{jk}(\theta_k - \theta_j)}\right] \label{ind}
\end{eqnarray}

To find the density $d\nu_n/d\lambda_N$, we use \eqref{nundef} and \eqref{ndef}:
\[
  \frac {d\nu_n} {d\lambda_N} = g(\tb) \cdot (2\pi)^N \cdot 
  \sum_{C \in \curs{V}} 
  \prod_{j,k=1}^N \frac{(t_j t_k)^{C_{jk} / 2} q_{jk}^{C_{jk}}}{C_{jk}!}
\]
To see that this is equal to \eqref{gffdens} and finish the proof, we use 
\eqref{ind}:
\begin{eqnarray*} 
  \lefteqn{ (2\pi)^N \sum_{C \in \curs{V}}
  \left( \prod_{j,k=1}^N \frac{(t_j t_k)^{C_{jk}/2}\, q_{jk}^{C_{jk}}}{C_{jk}!}
  \right) = \hspace{1.5in} } \\ \\ & = &
  \sum_{C \in \N^{N \times N} }
  \left[ \int_T d\btheta \prod_{l,m=1}^N e^{i C_{lm} (\theta_m - \theta_l)} 
  \right] \,
	\left[ \prod_{j,k=1}^N \frac{(\sqrt{t_j t_k} \, q_{jk})^{C_{jk}}}{C_{jk}!}
  \right]\\ & = &
  \int_T d\btheta \, \sum_{C \in \N^{N \times N}} 
	\left[ \prod_{j,k=1}^N e^{iC_{jk}(\theta_k - \theta_j)} \right] \,
  \left[ \prod_{j,k=1}^N\frac{(\sqrt{t_j t_k} \, q_{jk})^{C_{jk}}}{C_{jk}!}
  \right]\\ & = &
  \int_T d\btheta \Big[ \prod_{j,k=1}^N \sum_{C_{jk} \geq 0} \frac
  {\bigl( \sqrt{t_j t_k}q_{jk}\exp\{i(\theta_k - \theta_j)\} \bigr)^{C_{jk}}}
  {C_{jk}!} \Big] \\ & = &
	\int_T d\btheta \, \Big[\prod_{j,k=1}^N \exp \left\{ \sqrt{t_j t_k} \, 
	q_{jk}\,e^{i(\theta_k - \theta_j)} \right\} \Big].
\end{eqnarray*}

\end{document}